\documentclass[11pt,a4paper,final,reqno]{article}   % For Latex2e
\usepackage{amsthm, amsmath}   % For Latex2e
\usepackage{graphics, graphicx, amsfonts, epsf}   % For Latex2e
\theoremstyle{plain}
\newtheorem{theorem}{Theorem}[section]
\newtheorem{lemma}[theorem]{Lemma}
\newtheorem{proposition}[theorem]{Proposition}
\newtheorem{corollary}[theorem]{Corollary}
\newtheorem{observation}[theorem]{Observation}
\newtheorem{remark}[theorem]{Remark}

\newtheorem{example}[theorem]{Example}
\newtheorem{algorithm}[theorem]{Algorithm}

\theoremstyle{definition}
\newtheorem{definition}[theorem]{Definition}

\title{From graphs to tensegrity structures: \\
Geometric and symbolic approaches}
\author{
\addtocounter{footnote}{1} Miguel de Guzm\'an \footnote{In
memoriam.} \and \addtocounter{footnote}{-2} David Orden
\footnote{Departamento de Matem\'aticas, Universidad de Alcal\'a,
e-mail: david.orden@uah.es Research partially supported by grants
MEC MTM2005-08618-C02-02 and CAM S-0505/DPI/000235.}}
\date{}

\begin{document}
\maketitle

\begin{abstract}
A form-finding problem for tensegrity structures is studied; given
an abstract graph, we show an algorithm to provide a necessary
condition for it to be the underlying graph of a tensegrity in
$\mathbb{R}^d$ (typically $d=2,3$) with vertices in general
position. Furthermore, for a certain class of graphs our algorithm
allows to obtain necessary and sufficient conditions on the relative
position of the vertices in order to underlie a tensegrity, for what
we propose both a geometric and a symbolic approach.
\end{abstract}



\section{Introduction}
\label{section:Introduction}

In this paper we study an instance of the so-called
\emph{form-finding problems} for tensegrity structures. These have
brought special attention both among mathematicians and engineers
since the seminal works of Kenneth Snelson around $1948$
(see~\cite{snelson:web}). Roughly speaking, a form-finding problem
for a tensegrity structure asks to determine a geometric
configuration of points and straight edges in $\mathbb{R}^d$
(typically $d=2,3$) such that the whole structure is in a
self-tensional equilibrium. The word \emph{tensegrity} was coined
from \emph{tension} and \emph{integrity} by Buckminster Fuller,
deeply impressed by Snelson's work.

Apart from a purely mathematical interest~
\cite{connelly-whiteley:tensegrity,roth-whiteley:tensegrity},
understanding these structures has applications to architecture
and structural engineering~\cite{szabo-kollar:roofs} and has led
to interesting models for viruses and cellular structures
~\cite{caspar-klug:viruses,ingber:cellular}. It is also considered
a useful tool for the study of deployable structures
~\cite{motro:book,skelton:deployable,tibert:thesis}. Previous
works have proposed a number of different approaches to solve
form-finding problems, which can be found in the recent review
~\cite{pellegrino:review}.
%{\it OBS.: Buscar datos de esta referencia}.

In particular, the present paper deals with the form-finding
problem of building tensegrity structures with a given underlying
graph $G$, in a given $\mathbb{R}^d$. The graph has to be
understood as an \emph{abstract graph}, i.e., a set of vertices
and pairs of vertices (edges). We aim to solve the following two
problems:

\begin{itemize}
\item First, to decide whether $G$ can be the underlying graph of a
tensegrity structure in~$\mathbb{R}^d$.
\item In case such a tensegrity with underlying graph $G$ is possible,
to characterize the relative position of its vertices.
\end{itemize}

In order to solve these problems, we first look for decompositions
of tensegrities into basic instances, called \emph{atoms}. This
motivates a combinatorial method that allows to decompose a graph
$G$ into the smallest graphs that can underlie a tensegrity. In
order to build up a tensegrity with graph $G$, we propose to reverse
its decomposition: We show that this solves the above problems for a
certain class of graphs and we present two different approaches. The
first one looks at the geometric structure of the tensegrity; it is
quite visual and provides intuition of the intrinsic properties of
tensegrity structures. However, it becomes difficult to use for
complicated structures. The second approach condenses in a matrix
the information about the tensegrity; this allows to use tools from
Symbolic Computation, despite being less intuitive.

The paper is organized as follows: The basic notions and results are
introduced in Section~\ref{section:Preliminaries}. Then,
Section~\ref{section:Geometric approach} introduces a method to
decompose a tensegrity into atoms, which motivates a decomposition
of the abstract graph $G$, reversed then by geometric means.
Finally, in Section~\ref{section:Symbolic approach} a \emph{rigidity
matrix} is used for a symbolic resolution.
%
%In order to answer the first question, we present a combinatorial
%algorithm which decomposes the graph into the smallest possible
%tensegrities. The reversal of this decomposition leads to
%conditions on the vertices, for what we present two different
%approaches: The first method is geometric and hence more visual,
%but also more tricky. A symbolic method is proposed as a second
%choice; it is not so intuitive but puts together all the
%calculations in a matricial form and takes advantage of tools from
%Symbolic Computation.
%
%The reversal of each step of the algorithm consists in looking for
%the feasible positions of a new point with respect to the already
%placed ones. In order to do so, we propose two different
%approaches: On the one hand, the necessary and sufficient
%conditions can sometimes be deduced from the geometric structure
%of the tensegrity by taking the new point to limit positions.
%
%This is more visual; in essence, a tensegrity with the given
%underlying graph is constructed piece by piece and even a physical
%model can be built. But one has to handle a number of calculations
%in order to prove the experimental facts observed. On the other
%hand, the second method is not so visual, but uses a compact
%matrix formulation in order to perform the calculations needed for
%a proof. We show how to take advantage of Symbolic Computations
%using {\tt Maple}~\cite{maple} software.
%%{\it OBS.: comprobar referencia}.

\section{Preliminaries}
\label{section:Preliminaries}

In this section we introduce the basic notions and results used in
the paper. Despite it aims to be self-contained, an interested
reader can look at~\cite{whiteley:handbookDCG} for further
examples and a more detailed overview of the mathematical
concepts. Let us introduce first the rigorous definition of
``self-tensional equilibrium":

\begin{definition}\label{def:self-stress}
Let $G=(V,E)$ be an abstract graph:
\begin{itemize}
\item A \emph{framework} $G(P)$ in $\mathbb{R}^d$ is an embedding of
$G$ on a finite point configuration
$P:=\{p_1,\ldots,p_n\}$ in $\mathbb{R}^d$, with straight edges. In
the sequel we will focus on \emph{general position} point
configurations (no $d+1$ points lie on the same hyperplane).
% Edges $p_ip_j$ will also be denoted $ij$.
\item A \emph{stress} $w$ on a framework is an assignment of scalars
$w_{ij}$ (called \emph{tensions}) to its edges. Observe that
$w_{ij}=w_{ji}$, since they refer to the same edge.
\item Such a $w$ is called a \emph{self-stress} if, in addition, the following
equilibrium condition is fulfilled at every vertex:
\begin{equation}
\label{eq:equilibrium}
\forall i,\quad\sum_{ij \mbox{ edge}}w_{ij}(p_i-p_j)=0
\end{equation}
That is, for each vertex $p_i$ the scaled sum of incident vectors
$\stackrel{\longrightarrow}{p_ip_j}$ is zero.
\end{itemize}
\end{definition}

Observe that the null stress is always a self-stress, of no interest
for us. Note also that all scalar multiples of a self-stress (in
particular its opposite) are self-stresses as well. We will see
later that, indeed, the space of self-stresses on a given graph is a
vector space.

\begin{lemma}\label{lemma:incidence.degree}
Let $p\in P$ be a vertex of a $d$-dimensional framework $G(P)$
such that $P$ is in general position. Given a non-null self-stress
on $G(P)$, either at least $d+1$ of the edges incident to $p$
receive non-null tension, or all of them have null tension.
\end{lemma}

\begin{proof}
The result is true for any $d$, but the reader may consider
$d=2,3$ here. Let $k$ be the number of edges incident to $p$ that
have non-null tension. The equilibrium condition on $p$ implies
having $k$ vectors in $\mathbb{R}^d$, with common tail, which add
up to the zero vector. For $k<d+1$, this is only possible if their
$k+1$ endpoints do not span a $k$-space. But either $k=0$ or this
contradicts the general position assumption.
\end{proof}

As a consequence, the next property makes particularly interesting
the study of a certain family of general position frameworks, the
so-called \emph{$(d+1)$-regular} ones, for which a null tension on
a single edge propagates to the rest of them:

\begin{corollary}\label{coro:regular.non-null}
Given a $d$-dimensional framework all of whose points are in
general position and have exactly $d+1$ incident edges, a
self-stress is non-null if, and only if, it is non-null on every
edge.
\end{corollary}

The following definition introduces our final object of study,
which is a physical model of the mathematical objects defined
above:

\begin{definition}\label{def:tensegrity.structure}
We define a \emph{tensegrity structure} $T(P)$ to be a
self-stressed framework in which:
\begin{itemize}
\item Edges $ij$ such that $w_{ij}>0$ have been replaced by
inextensible \emph{cables} (its endpoints constrained not to get
further apart),
\item Edges with $w_{ij}<0$ have been replaced by unshrinkable
\emph{struts} (endpoints constrained not to get closer together),
and
\item Edges with $w_{ij}=0$ have been removed.
\end{itemize}
%\begin{itemize}
%\item Inextensible \emph{cables} (its endpoints
%constrained not to get further apart).
%\item Rigid \emph{bars} (constraining the distance between endpoints to stay the same).
%\item Unshrinkable \emph{struts} (endpoints constrained not to get
%closer together).
%\end{itemize}
\end{definition}

\begin{figure}[htb]
\begin{center}
%\vspace{-2.5cm}
%\epsfxsize=3in\leavevmode\epsfbox{signs.eps}
\includegraphics*[scale=0.4]{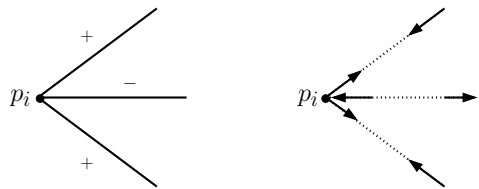}
%\vspace{-1.5cm}
\end{center}
\caption{ \label{fig:signs} Left: two cables ($+$) and a strut ($-$)
incident to $p_i$. Right: Their representation as springs with
inwards and outwards tensions.}
\end{figure}

If no confusion is possible, a tensegrity structure $T(P)$ will be
denoted by just $T$. For another physical interpretation, one can
think of cables and struts as springs endowed with a certain
tension, respectively inwards and outwards. That is; cables and
struts incident to point $p_i$ have respectively tensions in the
direction of {\tiny +\!\!} $\stackrel{\longrightarrow}{p_ip_j}$
(inwards) and -$\stackrel{\longrightarrow}{p_ip_j}$ (outwards),
see Figure~\ref{fig:signs}.

Observe that, given a tensegrity structure, it might be possible
to replace the struts by bars which react to the surrounding
tensions. For example, if we replace the strut in
Figure~\ref{fig:signs} by a bar, this will receive an outwards
tension at $p_i$, as a reaction to the sum of cable tensions. Such
a replacement is usual when constructing tensegrity sculptures,
like those in~\cite{snelson:web}.

\begin{figure}[htb]
\begin{center}
%\vspace{-2cm}
\includegraphics*[scale=0.5]{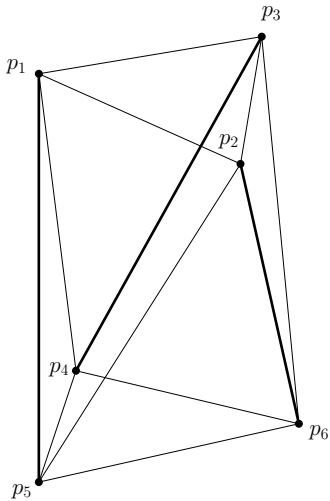}
%\epsfxsize=1.5in\leavevmode\epsfbox{otprism.eps}
%\vspace{-0.25cm}
\end{center}
\caption{ \label{fig:otprism} The oblique triangular prism.}
\end{figure}

The most emblematic tensegrity structure is shown in
Figure~\ref{fig:otprism}. Named \emph{oblique triangular prism}
with rotational symmetry, it is composed of nine cables, six of
which form two copies of an equilateral triangle, the top one
rotated $30$ degrees, joined by three struts alternating the rest
of cables. Thick edges denote the struts, which could be replaced
by bars as before.

%Provided the above definitions, in this paper we deal with the
%\emph{form-finding problem} of building tensegrity structures with
%a given underlying graph.
%finding out which configurations of points, cables and struts (or
%bars) give rise to similar tensegrity structures.
%
%In particular, in this paper we deal with two examples in which we
%are given a $4$-regular framework in general position in
%$\mathbb{R}^3$, one of whose points has unknown coordinates. The
%task is to find those positions of the undetermined point for
%which a non-null self-stress is admitted and hence a tensegrity
%structure with that topological configuration can be built. For
%the geometric treatment of more types of form-finding problems we
%refer the reader to~\cite{guzman:tensegrity}. Its algebraic
%resolution is postponed for a sequel paper.
%
%That is, when selected distance constraints on a finite number of
%points fix the point configuration (up to trivial motions in
%$\mathbb{R}^d$).
%Again , one can also ask for the topology or the combinatorics of
%the point configuration to be fixed, at least for small
%perturbations, see \cite{connelly-whiteley:tensegrity}, but in this
%paper we will focus on the first problem.
%
The last definition in this section introduces the smallest
tensegrities possible, which we will show in
Section~\ref{section:Geometric approach} to be the fundamental
bricks for building up any tensegrity:
%frameworks in general position that can admit a non-null
%self-stress:

\begin{definition}\label{def:atoms}
We define a \emph{self-stressed atom} in $\mathbb{R}^d$ ($d=2,3$) to
be a general position realization of the complete graph $K_{d+2}$,
together with its unique (up to constant multiplication) non-null
self-stress. The \emph{tensegrity atoms} are then obtained replacing
edges by cables and struts. When no confusion is possible, we will
just refer to \emph{atoms}.
\end{definition}

%\begin{figure}[htb]
%\begin{center}
%         \epsfxsize=5 in\leavevmode \epsfbox{figur1.eps}
%\end{center}
%\caption{How to obtain $\tilde B$ from $B$.} \label{fig.bloques}
%\end{figure}
\begin{figure}[htb]
\begin{center}
%\vspace{-2cm}
%\epsfxsize=2.1in\leavevmode\epsfbox{tensgratoms.eps}
\includegraphics*[scale=0.49]{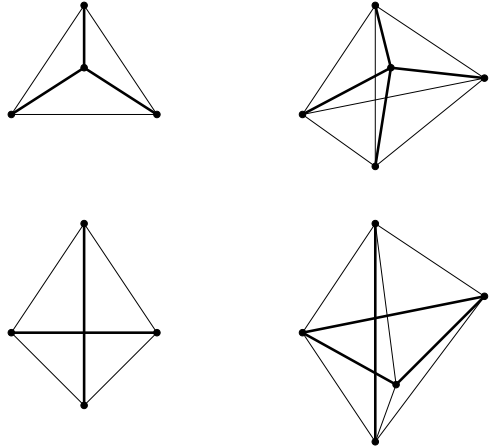}
%\vspace{-1.5cm}
\end{center}
\caption{ \label{fig:tensgratoms} Types of tensegrity atoms in
$\mathbb{R}^2$ (left) and $\mathbb{R}^3$ (right).}
\end{figure}

Figure~\ref{fig:tensgratoms} shows half of the possible tensegrity
atoms in $\mathbb{R}^d$ for $d=2,3$, where thick edges denote
struts. The other half is obtained by interchanging cables and
struts or, equivalently, by considering the opposite tensions. It is
not difficult to check, using Lemma~\ref{lemma:incidence.degree},
that configurations with fewer points or edges do not admit
self-stresses apart from the null one. The non-trivial fact that the
above frameworks do admit a unique (up to constants) non-null
self-stress appears in~\cite{rote-santos-streinu:ppt-polytope},
where existence is proved by the following result.

\begin{proposition}
\label{prop:dependence.stress} Let $\sum_{i=1}^{n} \lambda_i
p_i=0$, $\sum_{i=1}^{n} \lambda_i=0$ be an affine dependence on a
point set $P=\{p_1,\ldots,p_n\}$. Then,
$w_{ij}:=\lambda_i\lambda_j$ defines a self-stress on the complete
graph $K(P)$.
\end{proposition}

\begin{proof}
For any $p_i\in P$, we have:
\[
%\begin{array}{c}
\sum_{ij\in K}
w_{ij}(p_i-p_j)=\sum_{j=1}^{n}\lambda_i\lambda_j(p_i-p_j)=\lambda_i
p_i\sum_{j=1}^{n}\lambda_j-\lambda_i\sum_{j=1}^{n}\lambda_j p_j,
%=0.
%\end{array}
\]
which equals zero.
\end{proof}

In order to prove the uniqueness up to constants, consider two
different self-stresses, one not a scalar multiple of the other.
Then some linear combination of them would cancel the tension at a
particular edge but not at all of them, in contradiction with
Corollary~\ref{coro:regular.non-null}. Note that we are using the
claimed fact that self-stresses form a vector space, as will be
shown at the beginning of the next subsection.

%Note that the fact that we admit an atom self-stress to be
%multiplied by a constant allows the tension of one atom edge to be
%fixed, and the others consequently inherated.

\section{Geometric approach}
\label{section:Geometric approach}

In this section we present a geometric algorithm to decompose a
tensegrity into atoms. This decomposition motivates a combinatorial
one, which opens the way towards the resolution of the two problems
posed in Section~\ref{section:Introduction}.

%In particular, it solves the first one and in the second part of the
%section we show that, under certain conditions, its geometric
%reversal allows to solve the second problem.

\subsection{Decomposing into atoms}
\label{subsection:Decomposing into atoms}

Let $G(P)$ and $G'(P')$ be two self-stressed frameworks such that
$P\cup P'$ is a point configuration in general position. Let $w$
and $w'$ be their self-stresses. For the framework $G(P)\cup
G'(P')$ obtained by union of vertices and edges, one can define
the sum of self-stresses $w+w'$ in the natural way: Assign tension
$w_{ij}+w'_{ij}$ to common edges $ij$ and maintain the initial
tension at the others.

It is easy to observe that equations~(\ref{eq:equilibrium}) are
fulfilled and hence $w+w'$ is indeed a self-stress. Furthermore, the
space of self-stresses on a given graph $G=(V,E)$ together with this
sum and the product by a scalar form a vector subspace of
$\mathbb{R}^E$, when the latter is identified with the space of all
self-stresses.

Abusing notation, we denote by $G+G'$ the self-stressed framework
obtained. Observe that, after this addition is performed, one can
appropriately replace edges by cables and struts in order to obtain
a tensegrity structure $T+T'$. Hence, the sum of tensegrities yields
another tensegrity.

\begin{observation}\label{obs:sum.tensegrity}
We will only consider this kind of addition when $P$ and $P'$ have
at least~$d$ points in common; otherwise we obtain either two
separate tensegrity structures or one of them hanging from the
other.
\end{observation}

The main result in this section states that, reciprocally, under
our conditions every tensegrity can be decomposed into a sum of
tensegrity atoms:

\begin{theorem}{\bf (Atomic decomposition of tensegrities)}
\label{thm:decomposition} Every non-null tensegrity structure
$T(P)$, $P$ in general position, is a finite sum of tensegrity
atoms. This decomposition is not unique in general.
\end{theorem}

\begin{proof}
Let $G(P)$ and $w$ be the framework and non-null self-stress
associated to $T(P)$. We show how to obtain, by addition of atoms, a
chain of non-null self-stresses $w'$ on $G(P)$ in which the number
of vertices with only null incident tensions (\emph{null} vertices)
is increased at each step. At the end we come up with a
self-stressed framework with only null vertices, so that the
original tensegrity $T$ will be the sum of the opposites of those
atoms that have appeared in the process.

Let us focus on the two-dimensional case, since the $d$-dimensional
one is carried out analogously: At each step, an arbitrary non-null
vertex $a\in P$ is chosen to be converted in a null one. By
Lemma~\ref{lemma:incidence.degree}, only the following two cases are
possible (see Figures~\ref{fig:Step1Decomp}
and~\ref{fig:Step2Decomp}):

\begin{itemize}
\item {\bf Type 1:} If exactly three incident edges $ab,ac,ad$ have
non-null tension, we consider the atom $K$ of vertices $a,b,c,d$.
Since this atom has a non-null self-stress $w^{K}$ which is unique
up to constants, we can choose $w^{K}_{ab}$ to be the opposite of
the tension assigned to edge $ab$ at the current stress $w'$, i.e.
$w^{K}_{ab}:=-w'_{ab}$. Because of the equilibrium at~$a$, it turns
out that also $w^{K}_{ac}=-w'_{ac}$ and $w^{K}_{ad}=-w'_{ad}$.
Therefore, adding $w^{K}$ to the current self-stress makes vertex
$a$ have only null tensions at incident edges (i.e. makes it
disappear from the induced tensegrity). See
Figure~\ref{fig:Step1Decomp}, where dashed interior edges in the
second picture are opposite to those in the first one. Note that at
$b,c,d$ the edges $bc$, $bd$ and $cd$ may have appeared with
non-null tension, but these extra edges do not affect~$a$. However,
we will be concerned about them later.
\end{itemize}

\begin{figure}[htb]
\begin{center}
%\vspace{-2.5cm}
%\epsfxsize=4in\leavevmode\epsfbox{Step1Decomp.eps}
\includegraphics*[scale=0.49]{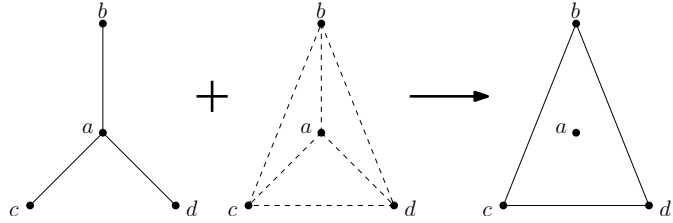}
%\vspace{-1.5cm}
\end{center}
\caption{ \label{fig:Step1Decomp} Type 1 step, exactly three
incident edges with non-null tension.}
\end{figure}

\begin{itemize}
\item {\bf Type 2:} If $a\in P$ has incidence degree greater than $3$,
let $b,c,d$ be neighbors of $a$. Consider the atom $\bar{K}$ of
vertices $a,b,c,d$ (and all the possible edges between them) and
choose it to have tension $w^{\bar{K}}_{ab}:=-w'_{ab}$ at edge $ab$.
Hence, obviously $w'+w^{\bar{K}}$ has null tension at edge $ab$.
Again, other edges $bc,bd,cd$ may appear with non-null tension, but
not incident to $a$. Hence, repeating this process if needed, we
obtain a self-stress on $G(P)$ in which $a$ has only three incident
edges with non-null tension (i.e. in the induced tensegrity, $a$ has
only three incident edges). Now we are in the previous case. See
Figure~\ref{fig:Step2Decomp}, where now only dashed edge $ab$ is
guaranteed to be opposite to its filled counterpart.
\end{itemize}

\begin{figure}[htb]
\begin{center}
%\vspace{-2.5cm}
%\epsfxsize=4in\leavevmode\epsfbox{Step2Decomp.eps}
\includegraphics*[scale=0.49]{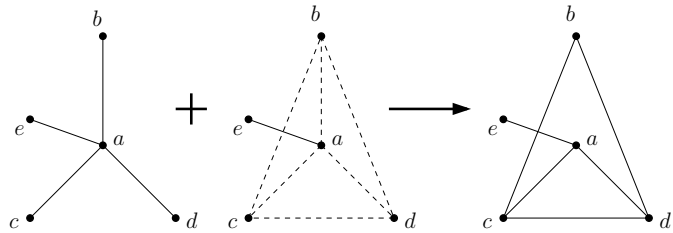}
%\vspace{-1.5cm}
\end{center}
\caption{ \label{fig:Step2Decomp} Type 2 step, more than three
incident edges with non-null tension.}
\end{figure}

Since a sum of self-stresses is another self-stress, after a finite
number of these steps we get a self-stress with at least one more
null vertex, for which we can iterate the process until all vertices
become null. Note that different choices of vertices to make null
may lead to different decompositions.
\end{proof}

%The reader should notice that adding an atom might cancel other
%tensions than the intended ones. The following lemma shows a
%restriction on this possibility:
%
\begin{remark}
\label{remark:atom.cancel}
%The only tensions that can be canceled by the addition of an atom are %precisely those of the atom.
%The addition of an atom cannot cancel other tensions that those of the atom %edges.
%The addition of an atom cannot cancel tensions of edges external to the atom.
The reader should notice that the addition of an atom only changes
the value of the self-stress on edges of the atom. Hence, adding an
atom might cancel other tensions than the intended ones, but only at
edges contained in the atom.
%The addition of an atom can only cancel the tensions of the
%edges contained in it.
\end{remark}
%
%\begin{proof}
%The proof is again provided for the case of dimension 2, although
%its generalization to general dimension $d$ is straightforward. For
%the sake of a contradiction, consider $K$ an atom of vertices
%$a,b,c,d$ and $\overline{pq}$ an edge with $\{p,q\}\not\subset
%\{a,b,c,d\}$, such that the addition of $K$ cancels the tension
%$w_{pq}$.
%
%Since $\overline{pq}$ is neither an edge of $-K$, the addition of
%this atom keeps $w_{pq}=0$ and therefore we conclude that the edge
%$\overline{pq}$ already had null tension before adding $K$, because
%the consecutive addition of an atom $K$ and its opposite $-K$ gives
%back the original self-stress.
%\end{proof}
%
%For an illustration of this proof, the reader can consider the graph
%in Figure~\ref{fig:BadCombinatorialDecomp}, where adding first the
%atom containing $a$ and then its opposite would lead to a tensegrity
%without edge $pq$ and therefore different from the original.

Motivated by the \emph{geometric} process in
Theorem~\ref{thm:decomposition}, we define now the following
\emph{combinatorial} algorithm, that can be applied to any abstract
graph~$G$:

\begin{algorithm}\textbf{(Combinatorial decomposition)}\\
\label{algo:combinatorial.decomposition}

\noindent INPUT: abstract graph $G=(V,E)$ and dimension $d$.\\
OUTPUT: list $L$ of ``atoms", where each atom is a subset of $(d+2)$
elements of $V$.
\begin{enumerate}
\item Initialize $L=\emptyset$.
\item While $E$ is not empty, choose a vertex $a\in V$ with
minimum degree and:
\begin{enumerate}
\item[2.1] If $a$ has degree $\leq d$, remove its incident edges
from $E$.
\item[2.2] If $a$ has degree $d+1$, let $a_0,\ldots,a_d$ be its
neighbors. Remove the edges $aa_i$ from $E$. Add to $E$ all the
edges $a_ia_j$ that were not in $E$. Insert the atom
$\{a,a_0,\ldots,a_d\}$ to the list $L$.
\item[2.3] If $a$ has degree at least $d+2$ do the following until
it has degree $d+1$, then go to $2.2$: Choose $d+1$ neighbors
$a_0,\ldots,a_d$ of $a$. Remove the edge $aa_0$ from $E$ and insert
to $E$ all the edges $a_ia_j$ that were not in $E$. Insert the atom
$\{a,a_0,\ldots,a_d\}$ to the list $L$.
\end{enumerate}
\item Return $L$.
\end{enumerate}
\end{algorithm}

%At each step, we choose a vertex and delete its incident edges using
%an \emph{atomic deletion operation}, defined as one of the
%operations depicted in Figures~\ref{fig:Step1Decomp}
%and~\ref{fig:Step2Decomp}, followed by a \emph{pruning operation}
%consisting in the removal of all the edges incident to vertices of
%degree smaller than $d+1$ in an iterative way, i.e. until no edge
%can be deleted because of this reason.

\begin{figure}[htb]
\begin{center}
%\vspace{-2.5cm}
%\epsfxsize=4in\leavevmode\epsfbox{Step2Decomp.eps}
\includegraphics*[scale=0.55]{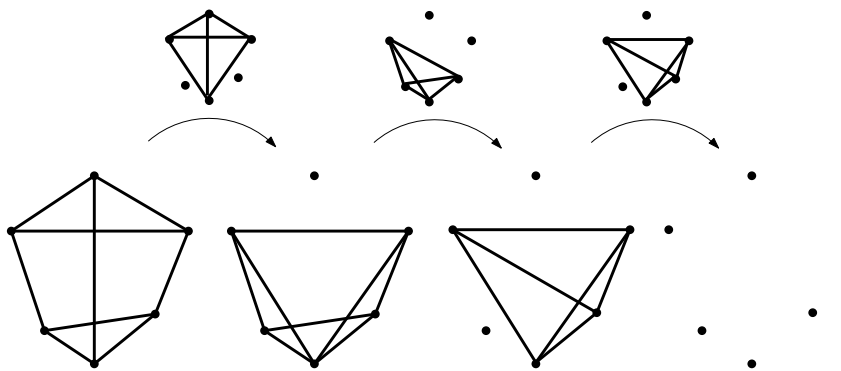}
%\vspace{-1.5cm}
\end{center}
\caption{ \label{fig:GoodCombinatorialDecomp} Combinatorial
decomposition fulfilling the condition of
Theorem~\ref{thm:algorithm-problem1}.}
\end{figure}

\begin{figure}[htb]
\begin{center}
%\vspace{-2.5cm}
%\epsfxsize=4in\leavevmode\epsfbox{Step2Decomp.eps}
\includegraphics*[scale=0.75]{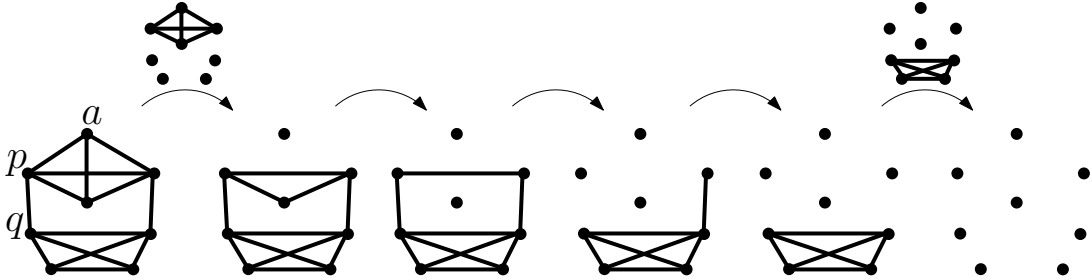}
%\vspace{-1.5cm}
\end{center}
\caption{ \label{fig:BadCombinatorialDecomp} Combinatorial
decomposition not fulfilling the condition of
Theorem~\ref{thm:algorithm-problem1}.}
\end{figure}

See Figures~\ref{fig:GoodCombinatorialDecomp}
and~\ref{fig:BadCombinatorialDecomp} for examples. This
combinatorial algorithm is the tool for the first result towards the
resolution of the problems posed in Section~1:

\begin{theorem} \label{thm:algorithm-problem1}
Given an abstract graph $G$, in order for it to underlie a
tensegrity $T(P)$ in~$\mathbb{R}^d$ it is a necessary condition
that, chosen a combinatorial decomposition, for every edge $pq$ of
$G$ there is an atom containing both endpoints $p$ and $q$. (In
other words, that step~$2.1$ of the algorithm removes only edges
that were contained in the complete graph defined by some atom of
the combinatorial decomposition).
\end{theorem}
\begin{proof}
Choose a combinatorial decomposition as above. If such a tensegrity
$T(P)$ exists, then the combinatorial decomposition induces a
geometric one: On the one hand, the geometric counterpart of step
$2.1$ shows up naturally because of
Lemma~\ref{lemma:incidence.degree}. On the other hand, steps $2.2$
and $2.3$ correspond to the geometric steps of types 1 and 2,
respectively, in the proof of Theorem~\ref{thm:decomposition}, with
tensions determined by the edge(s) to be deleted. The necessary
condition in the statement is then a consequence of
Remark~\ref{remark:atom.cancel}.
\end{proof}

Figure~\ref{fig:GoodCombinatorialDecomp} shows a combinatorial
decomposition that fulfills the condition of
Theorem~\ref{thm:algorithm-problem1}, while
Figure~\ref{fig:BadCombinatorialDecomp} shows one for which edge
$pq$ shows that the condition is not fulfilled, hence the graph
cannot underly a tensegrity.

Therefore, when asked about the existence of a tensegrity in
$\mathbb{R}^d$ with underlying graph $G$, we first get a
combinatorial decomposition and then check the condition of
Theorem~\ref{thm:algorithm-problem1}. If this is not fulfilled, then
the answer is negative. However, in the next subsection we show that
if the graph~$G$ admits a combinatorial decomposition of a certain
type, then this decomposition can be geometrically reversed, leading
in a straight way to a complete characterization of the point
sets~$P$ on which~$G(P)$ admits a tensegrity.

\subsection{Geometrically solving form-finding problems}
\label{subsection:Geometrically solving form-finding problems}

Let us start defining a class of combinatorial decompositions, for
which the reader can find an example in
Figure~\ref{fig:GoodCombinatorialDecomp}.

\begin{definition}
\label{def:edge-insert.comb.decomp} We call \emph{edge-inserting}
combinatorial decompositions those in which all extracted atoms but
the last one introduce at least one edge to the intermediate graph.
\end{definition}

The following result states that, for this class of combinatorial
decompositions, their geometric reversal provides a solution to the
problems posed in Section~\ref{section:Introduction}. Note the
request of a self-stress non-null on every edge, in order for $G$ to
underlie the tensegrity:

\begin{theorem}\label{thm:characterize}
If an abstract graph $G$ admits an edge-inserting combinatorial
decomposition, then the reconstruction of the graph from the atomic
decomposition produces a set of equations and negated equations
characterizing the point sets $P$ in general position that make
$G(P)$ admit a self-stress non-null on every edge (and hence the
ones that admit a tensegrity $T(P)$).
\end{theorem}
\begin{proof}
In order to determine which choices of coordinates underlie
tensegrities, we have to consider as variables the coordinates of
the new points $p_i$ added in the reconstruction (the first $d+2$
points can be arbitrarily chosen, since tensegrities are
projectively invariant~\cite{roth-whiteley:tensegrity}). Thus, the
tensions of the edges in each atom introduced are functions on these
variables.

Furthermore, recall from Definition~\ref{def:atoms} that the
self-stress of each atom has one degree of freedom. Hence, for the
reconstruction of a general combinatorial decomposition we have to
consider one extra variable $\alpha_j$ for each atom (except for the
first one, whose stress can be considered a normalization constant).

However, these extra variables are not needed for edge-inserting
combinatorial decompositions: For each step of the geometric
reconstruction, at least one edge was inserted by the combinatorial
decomposition algorithm and has to be removed in this precise
reconstruction step. Then, the relative self-stress given to this
atom (considered as a function on the positions) can be determined
at the insertion step: It is exactly the one that cancels the
tension(s) at the edge(s) to be removed.

Therefore, the reconstruction in this case leads to a system of
equations $f(p_1,\ldots,p_n)=0$ and negated equations
$f(p_1,\ldots,p_n)\neq 0$, one equation for each edge not in $G$
that appears in an intermediate step of the process (expressing that
this edge has tension zero and does not appear in the tensegrity)
and one negated equation for each edge in $G$ (expressing that it
does appear in the tensegrity).
\end{proof}

Let us point up that two decompositions of $G$ lead to equivalent
collections of conditions, since the tensegrities constructed using
one decomposition can always be decomposed and reconstructed using
the other. Hence, they have to fulfill both sets of necessary and
sufficient conditions, and therefore these have to be equivalent. In
order to finish the section, we illustrate the decomposition and
reconstruction process with the following example:

\begin{example}
\label{example:2d} The graph
$G=(\{1,\ldots,6\},\{12,14,16,23,26,34,35,45,56\})$ underlies a
tensegrity $T(P)$ in general position in $\mathbb{R}^2$ if, and only
if, the triangles $p_1p_2p_6$ and $p_3p_4p_5$ are in perspective
position, i.e., the lines $\overline{p_2p_3},\overline{p_5p_6}$ and
$\overline{p_1p_4}$ are concurrent.
\end{example}

\begin{figure}[htb]
\begin{center}
%\vspace{-2.5cm}
%\epsfxsize=3.7in\leavevmode\epsfbox{Example6points.eps}
\includegraphics*[scale=0.49]{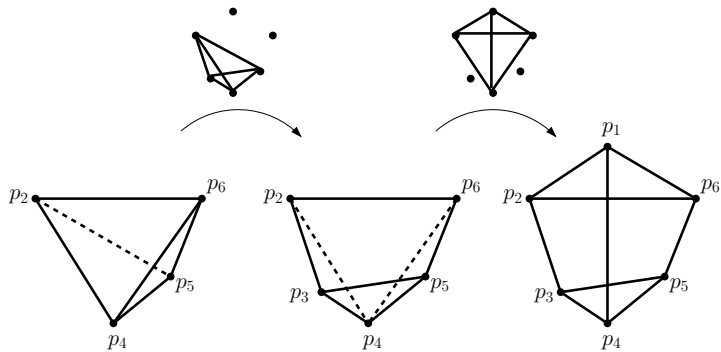}\\
%\vspace{-1.5cm}
\caption{ \label{fig:Example6points} Reconstruction of an
edge-inserting combinatorial decomposition.}
\end{center}
\end{figure}

We start with the edge-inserting combinatorial decomposition of $G$
depicted in Figure~\ref{fig:GoodCombinatorialDecomp}. In
Figure~\ref{fig:Example6points} we show the reversal of this
decomposition, where dashed edges are the inserted ones that have to
be removed.

As observed in the proof of Theorem~\ref{thm:characterize}, the
initial points $p_2,p_4,p_5,p_6$ can be chosen arbitrarily. Then we
add the atom $p_2,p_3,p_4,p_5$, whose self-stress is determined by
the condition of canceling the tension at edge $p_2p_5$. The
feasible positions of $p_3$ are characterized by the equations and
negated equations corresponding respectively to the edges omitted
and depicted in the middle picture.

Finally, the addition of the atom $p_1,p_2,p_4,p_6$ has to remove
edges $p_2p_4$ and $p_4p_6$. This fact allows us to rephrase
geometrically the set of equations and negated equations obtained
for the graph $G$: The point $p_1$ has to be placed on the line
spanned by the resultant of the ``undesired" tensions $w_{24}$ and
$w_{46}$ so that the tension on edge $p_1p_4$ replaces them. Moving
$p_1$ along the feasible line towards the intersection of
$\overline{p_2p_3}$ and $\overline{p_1p_4}$ forces tension $w_{26}$
at point $p_2$ to be very small. But this equals tension $w_{62}$ at
point~$p_6$; therefore, $p_1$ has to be close to the intersection of
$\overline{p_1p_4}$ and $\overline{p_5p_6}$ as well.

Taking this argument to the limit, we conclude that $p_1$ has to be
chosen in such a way that the lines
$\overline{p_2p_3},\overline{p_5p_6}$ and $\overline{p_1p_4}$ are
concurrent.

\section{Symbolic approach}
\label{section:Symbolic approach}

This section is devoted to show that the form-finding problems
considered can be reformulated in a matricial form, and how this can
be used to obtain the equations and negated equations in
Theorem~\ref{thm:characterize} by means of symbolic computations. In
particular, this provides an alternative to the geometric \emph{ad
hoc} reasonings like the one in Example~\ref{example:2d}, and allows
to deal with more complicated and less visual examples.

\subsection{The rigidity matrix}
\label{subsection:The rigidity matrix}

Before starting this subsection, let us point up that in this
paper we are not dealing with the rigidity of tensegrity
structures; we are just concerned by their self-equilibrium. An
interested reader can see~\cite{whiteley:handbookDCG} for more
information on rigidity of frameworks. However, the following
notion from rigidity analysis turns out to be useful for our
purposes:

\begin{definition}\label{def:rigiditymatrix}
Let $G(P)$ be a framework with $n$ vertices and $e$ edges in
$\mathbb{R}^d$. Its \emph{rigidity matrix} $R(P)$ has $e$ rows and
$nd$ columns, defined as follows:
\begin{itemize}
\item There is a row per edge $ij$ of the framework, with $i<j$ and in
lexicographic order.
\item Each block of $d$ columns is associated to a vertex $p_i$ and
contains either the $d$ coordinates $p_i-p_j$, at those rows
corresponding to edges $ij$ incident to $p_i$, or zeros at the
rest of rows.
\end{itemize}
For a \emph{complete framework} on $n$ vertices (with all the
possible edges), the rigidity matrix has the following condensed
form:
%\vspace{-0.5cm}
{\small
\[
\begin{array}{rll}
&
\begin{array}{cc}
\mbox{vertex $p_1$} &  \\
\downarrow &
\end{array}
&
\\

&
\left (
\begin{array}{ccccccc}
p_1\!-\!p_2 & p_2\!-\!p_1 &    0     &    0   & \ldots & 0 & 0 \\
p_1\!-\!p_3 &    0    & p_3\!-\!p_1  &    0   & \ldots & 0 & 0 \\
\vdots  & \vdots  & \vdots   & \vdots & \ddots & \vdots & \vdots\\
   0    &    0    &    0     &    0   & \ldots & p_{n-1}\!-\!p_n & p_n\!-\!p_{n-1} \\
\end{array}
\right )
&
\begin{array}{c}
\mbox{$\leftarrow$ edge $12$ }\\
\mbox{$\leftarrow$ edge $13$ }\\
\\
\\
\\
\end{array}
\end{array}
\]
}
\end{definition}

The key observation is that the equilibrium
equations~(\ref{eq:equilibrium}) can be restated in matricial form
as
\begin{equation}\label{eq:matricialequilibrium}
w\cdot R(P) = 0
\end{equation}
where $w$ is a $1\times e$ vector of entries $w_{ij},i<j$ (recall
that $w_{ij}=w_{ji}$) and the right-hand side is the $1\times nd$
zero vector. That is to say, self-stresses $w$ are row
dependencies for the rigidity matrix, what leads to the following
observation:

\begin{observation}
\label{obs:left.kernel} For a framework $G(P)$, being the underlying
graph of a tensegrity $T(P)$ is equivalent to the existence of a $w$
with no null component in the left kernel of the rigidity matrix
$R(P)$ of the framework.
\end{observation}

This observation already suggests a method to characterize the point
sets $P$ in general position that make $G(P)$ admit a tensegrity
$T(P)$: Computing the kernel of the rigidity matrix $R(P)$ for the
framework $G(P)$. For the sake of consistency, we refer the
interested reader to~\cite{guzman-orden-eaca04} for an example of
the use of this method, and we show here instead how to use the
decomposition-reconstruction method together with the matricial
expression of equilibrium.

\subsection{Symbolically solving form-finding problems}
\label{subsection:Symbolically solving form-finding problems}

%Our aim in this subsection is to obtain the equations and negated
%equations in Theorem~\ref{thm:characterize} by means of symbolic
%computations:
According to Observation~\ref{obs:sum.tensegrity}, the
reconstruction of a combinatorial decomposition inserts at most two
new points at each step, with indeterminate coordinates
$\overline{x},\overline{y}\in\mathbb{R}^d$. For the framework to
maintain equilibrium after their addition, the points have to belong
to the projection of the solutions
of~(\ref{eq:matricialequilibrium}) over the space of the
variables~$x_i,y_j$. Furthermore, in order to look for non-null
tensions on every edge, the inequations $w_{ij}\neq 0$ for every
$i,j$ have to be considered. Equivalently, extra variables $t_{ij}$
and equations $w_{ij}t_{ij}=1$ can be introduced to the system.

It is well-known (see Lemma~1 in page~120
of~\cite{cox-little-oshea:book}) that the polynomial elimination of
variables $w_{ij},t_{ij}$ from the above system of equations
contains the above projection. Therefore, polynomial elimination
provides necessary conditions for the position of the vertices,
which can be later symbolically tested for sufficiency.

We illustrate this symbolic method with the following example, in
which in order to perform polynomial elimination we have chosen to
compute the intersection of $\mathbb{R}[\overline{x},\overline{y}]$
with a \emph{Groebner basis} of the equations, for an elimination
order in which variables $x_i,y_j$ are smaller than the
$w_{ij},t_{ij}$'s; see Theorem~2 in page~113
of~\cite{cox-little-oshea:book}. Although there are more efficient
elimination techniques, the reason for this choice is to present the
operations performed with the software {\tt Maple}~\cite{maple}
which, apart from being broadly known, has a specific
\texttt{linalg} package for linear algebra. For more efficient
elimination software, the reader can consider {\tt
CoCoA}~\cite{cocoa}. Let us also point up that the graph considered
has eight triangles and, in fact, is the same as the graph of an
octahedron.

\begin{example}
\label{example:3d} The graph
$G=(\{1,\ldots,6\},\{12,13,14,15,23,25,26,34,36,45,46,56\})$ of an
oblique triangular prism underlies a tensegrity if, and only if, the
six vertices lie on a ruled hyperboloid that contains the edges of
one of the three cycles of length~$4$ in $G$. Equivalently, if and
only if the planes containing four alternating triangles intersect.
%four alternate faces of the oblique triangular prism
%intersect.
\end{example}

We start with the combinatorial decomposition of $G$ in which the
vertex~6 is chosen at step~2 of
Algorithm~\ref{algo:combinatorial.decomposition}. This gives the
atoms $2,3,4,5,6$ and $1,2,3,4,5$ (see Figure~\ref{fig:otprism}). In
the sequel we detail the steps of the reconstruction process using
{\tt Maple 7}, and omitting some outputs of no interest:

\bigskip

\noindent$\bullet$ According to the proof of
Theorem~\ref{thm:characterize}, the first $d+2$ points
$p_1,\ldots,p_5$ can be arbitrarily chosen:

\smallskip

\noindent{\tt >
p1:=[0,0,0]: p2:=[1,1,1]: p3:=[0,1,0]: %\\
p4:=[1,0,0]: p5:=[0,0,1]: %p6:=[x,y,z]:
}

\bigskip

\noindent$\bullet$ Solving the corresponding
equation~(\ref{eq:matricialequilibrium}) we get the tensions of the
atom $p_1,p_2,p_3,p_4,p_5$. In order to generate the equations, we
use the command {\tt geneqns} of the {\tt linalg} package, which
needs the equivalent transpose form $R^{t}\cdot w^{t}=0$
of~(\ref{eq:matricialequilibrium}). Then we solve them and, since
according to the proof of Theorem~\ref{thm:characterize} the stress
of the first atom can be considered a normalization constant, we
take the value 1 for the parameter obtained.
% We have chosen the names $w_{ij}$ for the unknown tensions:

\smallskip

%\noindent{\tt d:=3: n:=6: e:=12: zeros:=[0,0,0]: }
\noindent{\tt > d:=3: zeros:=[0,0,0]: with(linalg):}

\smallskip
\noindent{\tt > %\scriptsize
R:=matrix(10,5*d,[ \\
op(p1-p2),op(p2-p1),op(zeros),op(zeros),op(zeros),\\
op(p1-p3),op(zeros),op(p3-p1),op(zeros),op(zeros),\\
op(p1-p4),op(zeros),op(zeros),op(p4-p1),op(zeros),\\
op(p1-p5),op(zeros),op(zeros),op(zeros),op(p5-p1),\\
op(zeros),op(p2-p3),op(p3-p2),op(zeros),op(zeros),\\
op(zeros),op(p2-p4),op(zeros),op(p4-p2),op(zeros),\\
op(zeros),op(p2-p5),op(zeros),op(zeros),op(p5-p2),\\
op(zeros),op(zeros),op(p3-p4),op(p4-p3),op(zeros),\\
op(zeros),op(zeros),op(p3-p5),op(zeros),op(p5-p3),\\
op(zeros),op(zeros),op(zeros),op(p4-p5),op(p5-p4) ]): }

\smallskip
\noindent{\tt > Rt:=transpose(R):}

\smallskip
\noindent{\tt >
eqs:=geneqns(Rt,%\\
\verb+[w12,w13,w14,w15,w23,w25,w26,w34,w36,w45,w46,w56],+\\
vector(n*d,0)): }

\smallskip
\noindent{\tt >
tensions:=solve(eqs,{w12,w13,w14,w15,w23,w24,w25,w34,w35,w45}); }

\noindent $\mathit{tensions} := \{ w15 = -2\,w45, w35 = w45, w12 =
2\,w45, w14 = -2\,w45, w23 = -w45, \\ w34 = w45, w13 = -2\,w45, w45
= w45, w25 = -w45, w24 = -w45 \}$

\smallskip
\noindent{\tt > tensions:=subs(w45=1,tensions); }

\noindent $\mathit{tensions} := \{1 = 1, w35 = 1, w12 = 2, w14 = -2,
w23 = -1, w34 = 1, w13 = -2, \\
w25 = -1, w24 = -1, w15 = -2\}$

\bigskip

\noindent$\bullet$ Then we have to add the atom
$p_2,p_3,p_4,p_5,p_6$, in which $p_6$ has unknown
coordinates~$x,y,z$ and the stress of the atom is determined by the
cancelation of tensions at edges $\overline{p_2p_4}$ and
$\overline{p_3p_5}$. The same operations as above lead to a second
system of equations \textit{eqs2}, in which we substitute
$w24=1,w35=-1$ to get

\smallskip
\noindent $\mathit{eqs2subs} := \{ w23+w25-w26\,x+w26 = 0,
-w23-w34-x\,w36 = 0, w34+w45-w46\,x+w46 = 0, -w25-w45-x\,w56 = 0,
-w26+w26\,x+x\,w36-w46+w46\,x+x\,w56 = 0,
-w26+w26\,y-w36+w36\,y+y\,w46+y\,w56 = 0,
-w26+w26\,z+z\,w36+z\,w46-w56+w56\,z = 0, 1+w25-w26\,y+w26 = 0,
w23+1-w26\,z+w26 = 0, w34-1-w36\,y+w36 = 0, -w23+1-z\,w36 = 0,
-1-w34-y\,w46 = 0, -1-w45-z\,w46 = 0, -w25+1-y\,w56 = 0,
-1+w45-w56\,z+w56 = 0 \} $

\bigskip

\noindent$\bullet$ In order to look for non-null tension on every
edge, Corollary~\ref{coro:regular.non-null} turns out to be crucial;
if there is a non-null self-stress, then all its tensions are null.
%Therefore, the self-stress of the second atom is non-null since
%$w24=1$, and
Therefore, the self-stress obtained when adding the two atoms is
non-null since $w12=2\neq 0$ and this tension is not affected by the
second atom. In particular, it is not needed to consider extra
variables $t_{ij}$: In order to obtain the polynomial elimination,
we compute a Groebner basis of the polynomials in
$\mathit{eqs2subs}$, which we call $\mathit{polys2subs}$, for an
elimination order in which variables $x,y,z$ are smaller than the
$w_{ij}$'s. The command {\tt gbasis} of the {\tt Groebner} package
is used:

\smallskip

{\tt \noindent with(Groebner): \\
G:=gbasis(polys2subs,lexdeg([w23,w25,w26,w34,w36,w45,w46,w56],[x,y,z]));
}

\noindent $ \mathit{G} := [x^2-x-z^2-y^2+y+z, \ldots$
%1-w45+w56*z-w56, w25-1+y*w56, x*w56+w45+w25, 1+w45+z*w46,
%1+w34+y*w46, -w34-w45+w46*x-w46, y*w45+x-z, x*w45+z*w45-w45+y-1,
%w23-1+z*w36, -w34+1+w36*y-w36, w23+w34+x*w36, z*w34+x-y,
%x*w34+y*w34-w34+z-1, -w23-1+w26*z-w26, -1-w25+w26*y-w26,
%-w23-w25+w26*x-w26, z*w25-w25-x-y+1, x*w25-y*w25-z, y*w23-x-w23+1-z,
%x*w23-z*w23-y, w46*w56+w56*w45+w46*w45+w56-w46,
%w36*w56-w56-w46-w36-w26, w56*w34+w46,
%w26*w56+w56*w25+w26*w25+w56-w26, w56*w23+w26,
%w36*w46+w46*w34+w36*w34-w46+w36, w46*w26+w56+w46+w36+w26,
%w46*w25-w56, w46*w23-w36, w36*w45+w46, w45*w34-w46-w45-w34-1,
%w26*w45-w56, w45*w25+w56+w45+w25-1, w45*w23+1,
%w36*w26+w36*w23+w26*w23+w36-w26, w36*w25+w26, w26*w34-w36,
%w34*w25+1, w23*w34+w36+w34+w23-1, -w26-w25+w23*w25-1-w23
\quad and $40$ polynomials more, involving $w_{ij}$'s.

\bigskip

We conclude that the hyperboloid $x^{2}-y^{2}-z^{2}-x+y+z=0$
contains the set of points $p_6:=(x,y,z)$ such that the framework
admits a self stress non-null on every edge. Therefore this is a
necessary condition. In order to test its sufficiency, we have
considered an extra edge in the framework and forced its tension to
be null, obtaining the following stress from the left-kernel of the
corresponding rigidity matrix:
%
%\[
%w:= \left (
%\begin{array}{c}
%\frac{2y+2z-1-2zy-z^2+x^2-y^2}{-1+y-z+x}\\
%-\frac{2y+2z-1-2zy-z^2+x^2-y^2}{-1+y-z+x}\\
%-\frac{2y+2z-1-2zy-z^2+x^2-y^2}{-1+y-z+x}\\
%-\frac{2y+2z-1-2zy-z^2+x^2-y^2}{-1+y-z+x}\\
%-\frac{zx-zy-x+x^2-y^2+y}{-1+y-z+x}\\
%-\frac{-zy-z^2+z+xy-x+x^2}{-1+y-z+x}\\
%\frac{y+z-1+x}{-1+y-z+x}\\
%\frac{y(-y-z+1+x)}{-1+y-z+x}\\
%\frac{z-1+x-y}{-1+y-z+x}\\
%\frac{z(-y-z+1+x)}{-1+y-z+x}\\
%-\frac{-y-z+1+x}{-1+y-z+x}\\
%1
%\end{array}
%\right )
%\]

\[
w:= \left (
\begin{array}{c}
y+z-1-x\\
-y+x-z+1\\
-y+x-z+1\\
-y+x-z+1\\
-y+x\\
x-z\\
1\\
-\frac{y(-y+x-z+1)}{y+z+x-1}\\
-\frac{-y+z+x-1}{y+z+x-1}\\
-\frac{z(-y+x-z+1)}{y+z+x-1}\\
\frac{-y+x-z+1}{y+z+x-1}\\
-\frac{y-z-1+x}{y+z+x-1}
\end{array}
\right )
\]

We omit the (easy) computations checking that:
\begin{itemize}
\item Under the condition $x^{2}-y^{2}-z^{2}-x+y+z=0$ of the
hyperboloid, this is indeed a self-stress.
\item The denominator is not null: We observe that points
$p_3,p_4$ and $p_5$ already lie on the plane $y+z+x-1=0$.
Therefore, if point $p_6$ lied on that plane, the configuration
would not be in general position.
\item Similar arguments prove that all the components of $w$ are
non-null.
\end{itemize}

In conclusion, points $p_6=(x,y,z)$ for which a tensegrity exists
are precisely those lying on $x^{2}-y^{2}-z^{2}-x+y+z=0$. We now
observe that this hyperboloid contains the edges of the
quadrilateral $p_2p_3p_4p_5$. Indeed, it is the only one passing
through the initial five points and containing those four edges,
since the space of all hyperboloids has dimension~9. In order to
conclude the first condition in the statement we just have to
observe that the combinatorial decompositions that at step~2 of
Algorithm~\ref{algo:combinatorial.decomposition} choose vertices~1
or~6, lead to the~$4$-cycle $2345$, those choosing~3 or~5 lead to
$1264$, and those choosing~2 or~4 to $1365$ (see
Figure~\ref{fig:otprism}).

For the equivalent condition in the statement, which appears
in~\cite{crapo}, easy computations check that the four planes
defined by $p_1p_2p_3$, $p_1p_4p_5$, $p_3p_4p_6$ and $p_2p_5p_6$
intersect precisely for $x^{2}-y^{2}-z^{2}-x+y+z=0$, and the same
condition is obtained for the other four alternate planes.
%
%\begin{lemma}
%\label{lemma:equivalence.quadric} Let $p_1,\ldots,p_6$ be the six
%points considered above. Then, the following conditions are
%equivalent:
%\begin{enumerate}
%\item There is a quadric containing the six points and the edges of
%the quadrilateral $p_2p_3p_4p_5$.
%\item There is a quadric containing the six points and the edges of
%the quadrilateral $p_1p_4p_6p_2$.
%\item There is a quadric containing the six points and the edges of
%the quadrilateral $p_1p_4p_6p_2$.
%\end{enumerate}
%\end{lemma}
%
%%%%%%%%%%%%%%%%%%%%%%%%%%%%%%%%%%%%%%%%%%%%%%%%%%%%%%%%%%%%%%%%%%%%%%%%%%%%%%%%%%%%%
%
Let us finally note that an interested reader can find
in~\cite{montes:Comprehensive} a different symbolic approach to the
resolution of this problem, that uses comprehensive Groebner basis.

\begin{remark}
We have to call the attention of the reader to the fact that, in
spite of dealing with more complicated problems than the geometric
approach, the usefulness of the algebraic one is also limited by the
size and type of the problem. For instance, if
Corollary~\ref{coro:regular.non-null} cannot be used and all the
inequations $w_{ij}\neq 0$ have to be considered in order to look
for non-null tension on every edge. This would introduce more
auxiliary variables~$t_{12},t_{13},\ldots,t_{56}$ and a polynomial
$(w_{12}t_{12}-1)(w_{13}t_{13}-1)\cdots(w_{56}t_{56}-1)$, making the
computations infeasible.

Finally, let us recall that the above symbolic computations follow
the decomposition-reconstruction method from
Section~\ref{section:Geometric approach}. We already noted that
following Observation~\ref{obs:left.kernel} it is also possible to
compute symbolically the kernel of the $12\times 18$ rigidity
matrix~$R(P)$ for the whole framework $G(P)$. Instead, with the
decomposition-reconstruction approach we compute the kernel of the
two $10\times 15$ submatrices corresponding to the two atoms. The
main advantage of the decomposition reconstruction method is related
to this fact: Since $w24$ and $w35$ were no longer variables and the
system had fewer equations, the computation of the Groebner basis
was easier. In addition, we did not need to introduce an extra
variable~$t_{ij}$.
\end{remark}

\section*{Acknowledgements}
\label{section:Acknowledgements}

Both authors would like to thank the initial help and later
encouraging comments and revisions by professors Juan Llovet and
Tom\'as Recio, specially about the algebraic part. The second author
also wants to thank an anonymous reviewer for several valuable
comments that helped improving the results in the paper.

\section*{Final note}

The present paper is dedicated to the memory of Miguel de Guzm\'an,
who left us before it was finished. This work started when he gave a
conference about tensegrities at the Departamento de Matem\'aticas
of the Universidad de Alcal\'a in November~2003. Five months later
he died after a seminal version had been finished. Although part of
the results presented here still had to be formulated and proved in
their final form, their essence was already contained in that
preliminary version.


\begin{thebibliography}{10}

\bibitem{caspar-klug:viruses}
D.L.D.~Caspar and A.~Klug.
\newblock Physical principles in the construction of regular viruses.
\newblock In {\em Proceedings of Cold Spring Harbor Symposium on Quantitative Biology},
27:1--24, 1962.

\bibitem{cocoa}
CoCoATeam.
\newblock {\it CoCoA: a system for doing Computations in Conmutative
Algebra}.\\
{\tt http://cocoa.dima.unige.it}


\bibitem{connelly-whiteley:tensegrity}
R.~Connelly and W.~Whiteley.
\newblock Second-order rigidity and prestress stability for tensegrity frameworks.
\newblock {\em SIAM Journal of Discrete Mathematics}, 9(3):453--491, August 1996.

\bibitem{cox-little-oshea:book}
D.~Cox, J.~Little and D.~O'Shea.
\newblock {\em Ideals, varieties and algorithms.}
\newblock Springer Verlag, Berlin, 1991.

\bibitem{crapo}
H.~Crapo.
\newblock Structural rigidity.
\newblock {\em Structural Topology} 1:26--45, 1979.

\bibitem{guzman-orden-eaca04}
M.~de Guzm\'an and D.~Orden.
\newblock Finding tensegrity structures: Geometric and symbolic aproaches.
\newblock In {\em Proceedings of EACA-2004}, 167--172, 2004.

\bibitem{ingber:cellular}
D.E.~Ingber.
\newblock Cellular tensegrity: defining new rules of biological
design that govern the cytoskeleton.
\newblock {\em Journal of Cell Science}, 104:613--627, 1993.

\bibitem{montes:Comprehensive}
M.~Manubens and A.~Montes.
\newblock Improving DISPGB algorithm using the discriminant
ideal.
\newblock To be published at the special A3L issue of the
{\em Journal of Symbolic Computation}, 2006.
\newblock Available at \texttt{http://www.arxiv.org/abs/math.AC/0601763}

\bibitem{maple}
Waterloo Maple Inc.
\newblock {\it Maple}. {\tt http://www.maplesoft.com}

\bibitem{motro:book}
R.~Motro.
\newblock {\em Tensegrity: Structural systems for the future.}
\newblock Kogan Page Science, London, 2003.

\bibitem{pellegrino:review}
S.~Pellegrino and A.G.~Tibert.
\newblock Review of form-finding methods for tensegrity structures.
\newblock {\em International Journal of Space Structures}, 18(4):209--223(15), December 2003.

\bibitem{rote-santos-streinu:ppt-polytope}
G.~Rote, F.~Santos and I.~Streinu.
\newblock Expansive motions and the polytope of pointed pseudo-triangulations.
\newblock In B.~Aronov, S.~Basu, J.~Pach and M.~Sharir, editors, {\em Discrete and
  Computational Geometry -- The Goodman-Polack Festschrift}, volume~25 of
{\em Algorithms and Combinatorics}. Springer Verlag, Berlin, 2003,
699-736.

\bibitem{roth-whiteley:tensegrity}
B.~Roth and W.~Whiteley.
\newblock Tensegrity frameworks.
\newblock {\em Transactions of the American Mathematical Society}, 265:419--446, 1981.

\bibitem{skelton:deployable}
R.E.~Skelton.
\newblock Deployable tendon-controlled structure.
\newblock {\em United States Patent $5\,642\,590$}, July 1, 1997.

\bibitem{snelson:web}
K.~Snelson.
\newblock {\tt http://www.kennethsnelson.net}

\bibitem{szabo-kollar:roofs}
J.~Szabo and L.~Koll\'ar.
\newblock Structural design of cable-suspended roofs.
\newblock {\em Akademiai Kiado}, Budapest 1984.

\bibitem{tibert:thesis}
A.G.~Tibert.
\newblock Deployable tensegrity structures for space applications.
\newblock Ph.D. Thesis, {\em Royal Institute of Technology}, Stokholm 2002.

\bibitem{whiteley:handbookDCG}
W.~Whiteley.
\newblock Rigidity and scene analysis.
\newblock In J.~E. Goodman and J.~O'Rourke, editors, {\em Handbook of Discrete
  and Computational Geometry}, chapter~49, pages 893--916. CRC Press, New York,
  1997.

%\bibitem{salas:calculus}
%S.~Salas and E.~Hille.
%\newblock {\em Calculus: One and Several Variable}.
%\newblock John Wiley and Sons, New York, 1978.

\end{thebibliography}
\end{document}